\documentclass[11pt,draft]{article}
\usepackage{amsfonts}
\usepackage{mathrsfs}
\usepackage{amsmath,amsfonts,mathrsfs,amssymb,color}
\usepackage{indentfirst}

\numberwithin{equation}{section}

\newtheorem{theo.}{\quad\, Theorem}[section]
\newtheorem{defi.}{\quad\, Definition}[section]
\newtheorem{lemm.}{\quad\, Lemma}[section]
\newtheorem{coro.}{\quad\, Corollary}[section]

\begin{document}

\title{Monotone iterative technique for nonlinear fourth order integro-differential
equations$^*$ }
\author{Jinxiang Wang$^{\ast}$
\\
 {\small  Department of Applied Mathematics, Lanzhou University of Technology, Lanzhou, P.R. Chin}\\
}
\date{} \maketitle
\footnote[0]{E-mail address: wjx19860420@163.com(Jinxiang Wang),  
} \footnote[0] {$^*$Corresponding author: Jinxiang Wang.\  Supported by the NSFC(No.11961043,\ No. 11801453). }

 \begin{abstract}
\baselineskip 18pt
In this paper, we consider the solvability of a class of nonlinear fourth order integro-differential equations with Navier boundary condition.  We first deal with a corresponding linear problem and establish a maximum principle. Using the maximum principle, we develop a monotone iterative technique in the presence of lower and upper solutions to solve the nonlinear problem under certain conditions.  Some examples are presented to illustrate the main results.

 \end{abstract}


\vskip 3mm

{\small\bf MR(2000)\ \ \ 34B10, \ 34B18}

\baselineskip 22pt

\section{Introduction}

The aim of this paper is to develop a monotone iterative technique
\ for the following nonlinear fourth order Fredholm type Integro-differential equation(IDE)
$$y^{(4)}(x)=f(x,y(x),\int_0^1k(x,t)y(t)dt),\ \ \ \ x\in (0,1),
\eqno (1.1)
$$
with the Navier boundary condition
$$y(0) = y(1) = y''(0) = y''(1) = 0,
\eqno (1.2)
$$
where $f\in C([0,1]\times {\mathbb{R}}^{2}, \mathbb{R})$, $k\in C([0,1]\times [0,1], \mathbb{R})$.

The boundary value problem (1.1),(1.2) can be seen as a generalization of the linear fourth order problem
$$
\aligned
&y^{(4)}(x)+My(x)-N\int_0^1k(x,t)y(t)dt=p(x),\ \ \ \ x\in (0,1),\\
&y(0) = y(1) = y''(0) = y''(1) = 0\\
\endaligned
\eqno (1.3)
$$
where $M, N$ are constants, $p\in C[0,1]$. Problem (1.3) arises from  the models for suspension bridge [1-2], quantum
theory [3] and transient ultrasonic fields [4]. Since the nonlocal term under the integral sign will cause some mathematical difficulties, the analytical solutions for IDEs are usually not easy to obtain. 
For the linear fourth order boundary value problem governed by integro-differential equations like (1.3), only a few studies have been carried out by using numerical methods, see, e.g., [5-11] and  the references therein. 

The monotone iterative technique concerning upper and lower solutions is a powerful tool to solve the lower-order
differential equations with various kinds of boundary conditions, see e.g., [12-17] and the references therein. This technique is that, for the considered problem, starting from a pair ordered lower and upper, one constructs two monotone sequences such that them uniformly
converge to the extremal solutions between the lower and upper solutions. This technique has also been applied to the
special case of (1.1),(1.2) that $f$ does not contain integral term, namely simple fourth-order boundary value
problem
$$
\aligned
&y^{(4)}(x)=f(x,y(x)),\ \ \ \ x\in (0,1),\\
&y(0) = y(1) = y''(0) = y''(1) = 0,\\
\endaligned
\eqno (1.4)
$$
see [18-20].

It is worth noticing that, in [20], Cabada et al. pointed out that: for general second-order differential equation with periodic,
Neumann, or Dirichlet boundary conditions, it is well
known that the existence of a well-ordered pair of lower and upper solutions $\alpha\leq\beta$ is sufficient to ensure the existence of a solution in the sector enclosed by them. However, this result is not true for fourth-order differential equations (1.4), see the counterexample in [20, Remark 3.1]. Indeed, the application of the lower and upper solutions method in boundary value problems of fourth order is heavily dependent on the conclusion of the maximum principle for the corresponding linear operators. For fourth order local problems without integral terms similar to (1.4), other related results on the lower and upper solutions method and monotone iterative technique, see
[21-38] and the references therein.

As far as we know, there have been no studies on the lower and upper solutions method and monotone iterative technique
for nonlocal fourth order problems such as nonlinear IDE (1.1) with corresponding boundary value conditions. 
 Motivated by the above works described, the main object of this paper is to build maximum principle for the linear problem (1.3), and then develop a monotone iterative technique in the presence of lower and upper solutions to solve nonlinear problem (1.1),(1.2) under certain conditions.

The rest paper is arranged as follows: In Section 2, we first prove a uniqueness result of solutions for linear integro-differential equation (1.3) with inhomogeneous boundary value condition, and then we establish  maximum principle for the operator in (1.3).
In Section 3, 
based upon maximum principle, we develop a monotone iterative technique for (1.1),(1.2) in the presence of lower and upper solutions under some monotonic condition on the nonlinearity $f$. Finally, in Section 4, we present two examples to illustrate the main results. The first one is a concrete nonlinear nonlocal fourth order boundary value problem; The second one is a general sixth order boundary value problems which can be transformed into the fourth order nonlocal problem like (1.1),(1.2).

\vskip 3mm

\section{The linear inhomogeneous boundary value problem governed by integro-differential equations}

In this Section, we prove a uniqueness result of solutions for (1.3) with general inhomogeneous boundary value condition, and then establish maximum principle for the corresponding operators.

As preliminaries, we first consider the following linear fourth order inhomogeneous boundary value problem
 $$
\aligned
&y^{(4)}(x)+My(x)=p(x),\ \ \ \ x\in (0,1),\\
&y(0)=A, \ y(1)=B, \ y''(0)=C, \ y''(1)=D,\\
\endaligned
\eqno (2.1)
$$
where $M>0$ and $A, B, C, D$ are constants, $p\in C[0,1]$.
\ By Cabada [20, Lemma 2.1] or Ma et al. [23, Theorem 2.1], if $M\leq c_{0}\approx  950.8843$, then Problem (2.1) has a unique solution given by
$$y(x)=\int^{1}_{0}G(x,s)p(s)ds+Aw(x)+Bw(1-x)+C \chi(x)+D\chi(1-x), \ \ \ \  \ x\in [0,1],\eqno (2.2)
$$
where
$$G(x,s)=\left\{
\begin{aligned}\frac{1}{2m^{2}}
[\frac{\sin(mx)\sin(m(1-s))}{m\sin m}-\frac{\sinh(mx)\sinh(m(1-s))}{m\sinh m}],~~~0\leq x\leq s\leq1,\\
\frac{1}{2m^{2}}
[\frac{\sin(ms)\sin(m(1-x))}{m\sin m}-\frac{\sinh(ms)\sinh(m(1-x))}{m\sinh m}],~~~0\leq s\leq x\leq1,
\end{aligned}
\right.\eqno (2.3)$$
is the Green function for the boundary value problem $$
\aligned
&y^{(4)}(x)+My(x)=0,\ \ \ \ x\in (0,1),\\
&y(0) = y(1) = y''(0) = y''(1) = 0\\
\endaligned
\eqno (2.4)
$$
with $M=m^{4}$;
$w(x)$ is the unique solution of the inhomogeneous problem
$$y^{(4)}(x)+My(x)=0, \ \ \ \  \ \ y(0)=1, \ y''(0)=y(1)=y''(1)=0
\eqno (2.5)
$$
and $\chi(x)$ is the unique solution of the inhomogeneous
problem
$$y^{(4)}(x)+My(x)=0, \ \ \ \ \ \ y(0)=0, \ y''(0)=1, \ y(1)=y''(1)=0.
\eqno (2.6)
$$
Moreover, by Cabada [20, Proposition 2.1.], when $0<M=m^{4}\leq c_{0}\approx  950.8843$, the Green function $G(x,s)$ given by (2.3) is nonnegative on $[0,1]\times [0,1]$. 

Denote $$h(x,A,B,C,D)= Aw(x)+Bw(1-x)+C \chi(x)+D\chi(1-x),\eqno (2.7)$$ it is easy to see that $h$ is the unique solution of
 $$
 \aligned
 &y^{(4)}(x)+My(x)=0, \ \ \ \ x\in (0,1),\\
 &y(0)=A, \ y(1)=B, \ y''(0)=C, \ y''(1)=D,\\
 \endaligned
$$
and the following conclusion hold:
\vskip 3mm

\noindent{\bf Lemma 2.1}\ Assume that $0<M\leq c_{1}\approx  125.137 $. Then

if $A\geq0, B\geq0, C\leq 0, D\leq 0$, we have $$h(x,A,B,C,D)= Aw(x)+Bw(1-x)+C \chi(x)+D\chi(1-x)\geq0;\eqno (2.8)$$

if $A\leq0, B\leq0, C\geq 0, D\geq 0$, we have $$h(x,A,B,C,D)= Aw(x)+Bw(1-x)+C \chi(x)+D\chi(1-x)\leq0.\eqno (2.9)$$

\noindent{\bf Proof.}\ By Cabada [20, Corollary 2.1.], if $0<M\leq c_{1}\approx  125.137 $, then $w(x)\geq 0, \forall x\in [0,1]$ and $\chi(x)\leq 0, \forall x\in [0,1]$, and thus (2.8) and (2.9) are immediate consequences. \hfill{$\Box$}

\vskip 3mm
Now, we give out the first main result of this section.

\noindent{\bf Theorem 2.1}\ Assume that $0<M\leq c_{0}\approx  950.8843$, and $$\|k\|_{\infty}=\max\{k(x,t)|\ (x,t)\in [0,1]\times[0,1]\} < \frac{1}{|N| \underset{x\in [0,1]}\max\int^{1}_{0}G(x,s)ds},\eqno (2.10)$$ then for any $p\in C[0,1]$ and constants  $A, B, C, D$, the following nonlocal inhomogeneous boundary value problem
$$
\aligned
&y^{(4)}(x)+My(x)-N\int_0^1k(x,t)y(t)dt=p(x),\ \ \ \ x\in (0,1),\\
&y(0)=A, \ y(1)=B, \ y''(0)=C, \ y''(1)=D\\
\endaligned
\eqno (2.11)
$$
has a unique solution.

\noindent{\bf Proof.} Observe that $y$ is a solution of (2.11) if and only if $y$ is a fixed point of the the operator $K: C[0,1]\rightarrow C[0,1]$ given by
$$
[Ky](x)=\int^{1}_{0}G(x,s)p(s)ds +N\int^{1}_{0}\int^{1}_{0}G(x,s)k(s,t)y(t)dtds+h(x,A,B,C,D),
$$
where $G(x,s)$ and  $h$ is as in (2.3) and (2.7) respectively.

For $u,v\in C[0,1]$, by (2.10) we have
$$
\aligned
\|Ku-Kv\|_{\infty}=&\|N\int^{1}_{0}\int^{1}_{0}G(x,s)k(s,t)[u(t)-v(t)]dtds\|_{\infty}\ \\
&\leq \|u-v\|_{\infty}|N| |\int^{1}_{0}\int^{1}_{0}G(x,s)k(s,t)dtds|\\
&\leq \|u-v\|_{\infty}|N| \|k\|_{\infty}\int^{1}_{0}G(x,s)ds\\
&\leq \|u-v\|_{\infty}|N| \|k\|_{\infty}\underset{x\in [0,1]}\max\int^{1}_{0}G(x,s)ds\\
&< \|u-v\|_{\infty}.\\
 \endaligned
  $$
Based on Banach fixed point theorem, there exists a  unique fixed point for the operator $K$, which assures the existence and uniqueness of
solution for (2.11).  \hfill{$\Box$}

\vskip 3mm

\noindent{\bf Remark 2.1}  When $A=B=C=D=0$, the inhomogeneous problem (2.11) will degenerate into the homogeneous problem (1.3). Thus we obtain a result of the existence and uniqueness of solutions for (1.3). As far as we know, this result is new.

In the sequal, we establish maximum principle for the operator in (2.11) and (1.3).

Using Picard's iterative method we know that for any $y_{0}\in C[0,1]$, the sequence given by $y_{n}=Ky_{n-1},\ n\geq1$, converges to the unique solution given by Theorem 2.1. Taking
$$
y_{0}=\int^{1}_{0}G(x,s)p(s)ds+h(x,A,B,C,D)
$$
we get that
$$
y_{n}(x)=y_{0}(x)+\int^{1}_{0} Q_{n}(x,s)h(s,A,B,C,D)ds+\int^{1}_{0} F_{n}(x,s)p(s)ds,
$$
where
$$
F_{n}(x,s)=\int^{1}_{0}Q_{n}(x,t)G(t,s)dt,
$$
$$
Q_{n}(x,s)=\sum_{i=1}^{n}R^{(i)}(x,s).
$$
Here
$$
R^{(i)}(x,s)=\int^{1}_{0} R^{(i-1)}(x,t)R(t,s)dt,\ \ \ i\geq2,
$$
and
$$
R^{(1)}(x,s)=R(x,s)=N\int^{1}_{0}G(x,t)k(t,s)dt.
$$
By (2.10) we have $$N \int^{1}_{0}G(x,s)k(s,t)dt\leq|N| \|k\|_{\infty}\underset{x\in [0,1]}\max\int^{1}_{0}G(x,s)ds=d<1,$$
then
 $$
 \|R^{(i)}\|_{\infty}\leq d^{i}
 $$
 and the series $\sum_{i=1}^{\infty}R^{(i)}(x,s)$ will converge to a function $Q\in C([0,1]\times[0,1])$. Meanwhile, $\{F_{n}(x,s)\}$ will converge to the function $F\in C([0,1]\times[0,1])$ given by
 $$
 F(x,s)=\int^{1}_{0}Q(x,t)G(t,s)dt.
 $$
 Now, by passing  to the limit for the Picard's iterative $y_{n}=Ky_{n-1}$, we conclude that the unique solution of (2.11) is given by
$$
y(x)=\int_0^1G(x,s)p(s)ds+h(x,A,B,C,D)+\int^{1}_{0}Q(x,s)h(s,A,B,C,D)ds+\int_0^1F(x,s)p(s)ds. \eqno (2.12)
$$

\vskip 3mm
Now, by Lemma 2.1, the expression (2.12), and the positivity of the Green function $G$, one can easily get the following maximum principle for problem (2.11).
\vskip 3mm

\noindent{\bf Theorem 2.2}\  Assume that $0<M<c_{1}\approx  125.137$, $Nk(x,t)\geq0$ and $\|k\|_{\infty}< \frac{1}{|N| \underset{x\in [0,1]}\max\int^{1}_{0}G(x,s)ds}$, then

(i)\ \  If $p(x)\geq0, \forall x\in [0,1]$, $A\geq0, B\geq0, C\leq 0, D\leq 0$, the unique solution of (2.11) $y\geq0$;

(ii)\ \ If $p(x)\leq0, \forall x\in [0,1]$, $A\leq0, B\leq0, C\geq 0, D\geq 0$, the unique solution of  (2.11) $y\leq0$.

\vskip 3mm

\noindent{\bf Remark 2.2} The conclusion of Theorem 2.2 also hold for homogeneous problem (1.3). That is, we get a maximum principle for the fourth-order differential operator $Ly=y^{(4)}+My-N\int_0^1k(x,t)y(t)dt$ in function space $D(L):=\{y\in C^{4}[0,1]:y(0)=y(1)=y''(0)=y''(1)=0\}.$

\section{Main results}

\noindent{\bf Definition 3.1.} The function $\alpha\in C^4 [0, 1]$ is said to be a {\it lower solution} for the BVP (1.1),(1.2) if
$$\alpha^{(4)}(x)\leq f(x,\alpha(x),\int_0^1k(x,t)\alpha(t)dt),\ \ \ \ x\in (0,1),
\eqno (3.1)
$$
and
$$
\alpha(0)\leq 0,\ \alpha(1) \leq 0,\ \alpha''(0) \geq 0, \ \alpha''(1) \geq 0.
\eqno (3.2)
$$
An {\it upper solution} $\beta \in C^4[0, 1]$ is defined analogously by reversing the inequalities in (3.1),(3.2).

\noindent{\bf Theorem 3.1.} Assume $f\in C([0,1]\times {\mathbb{R}}^{2}, \mathbb{R})$, $k\geq0$ and there exist a lower solution $\alpha$ and an upper solution $\beta$ for problem (1.1), (1.2) which satisfy
$$\alpha(x) \leq \beta(x) \ \ \ \ \text{for}\ x \in [0, 1].
\eqno (3.3)
$$
If there exist two constants $M, N>0$ satisfying $M<c_{1}\approx  125.137$ and $\|k\|_{\infty}< \frac{1}{N \underset{x\in [0,1]}\max\int^{1}_{0}G(x,s)ds}$ such that
$$
f(x,u_{1},v_{1})-f(x,u_{2},v_{2})\geq -M(u_{1}-u_{2})+N(v_{1}-v_{2})\eqno (3.4)
$$
for $$\ \alpha(x) \leq u_1 \leq u_2 \leq \beta(x) \ \ \ \ \text{and}\ \ \ \int_0^1k(x,t)\alpha(t)dt \leq v_1 \leq v_2 \leq \int_0^1k(x,t)\beta(t)dt,\ \ \ \
  x\in [0,1].
$$
Then the iterative sequences $\{\alpha_{n}\}$ and $\{\beta_{n}\}$ produced by the iterative procedure
$$
\aligned
&Ly_{n}(x)=f(x,y_{n-1}(x),\int_0^1k(x,t)y_{n-1}(t)dt)+My_{n-1}(x)-N\int_0^1k(x,t)y_{n-1}(t)dt,\\
&y_{n}(0) = y_{n}(1) = y_{n}''(0) = y_{n}''(1) = 0,\ \ \ \ \ n=1,2,\ldots,\\
\endaligned
\eqno (3.5)
$$
with the initial functions $y_{0}=\alpha$ and $y_{0}=\beta$ respectively satisfy
$$
\alpha_{n-1} \leq \alpha_{n}\leq \beta_{n}\leq \beta_{n-1}\eqno (3.6)
$$
and converge uniformly to the extremal solutions of BVP (1.1),(1.2) in $[\alpha,\beta]$. 

\noindent{\bf Proof.} Define
the mapping $E: C[0,1]\rightarrow C[0,1]$  by
$$
E(\sigma)(x)=f(x,\sigma(x),\int_0^1k(x,t)\sigma(t)dt)+M\sigma(x)-N\int_0^1k(x,t)\sigma(t)dt
$$ Denote $\Phi=T\circ E$, where $T=L^{-1}: C[0,1]\rightarrow D(L)$. By Remark 2.1, it is easy to see that $T: C[0,1]\rightarrow C[0,1]$ is compact, then $\Phi: C[0,1]\rightarrow C[0,1]$ is completely continuous. Obviously, the solutions of (1.1),(1.2) in $C[0,1]$ is equivalent to the fixed-points of the mapping $\Phi$.

 Firstly, we show that
 $$\alpha\leq u\leq\beta\Rightarrow \alpha\leq \Phi(u)\leq\beta.\eqno (3.7)$$
 Let $g=\Phi u-\alpha$, by the definition of the lower solutions and (3.4), we have
$$
\aligned
&g^{(4)}+Mg-N\int_0^1k(x,t)g(t)dt\\
&=[(\Phi u)^{(4)}+M(\Phi u)-N\int_0^1k(x,t)(\Phi u)(t)dt]-[\alpha^{(4)}+M\alpha-N\int_0^1k(x,t)\alpha(t)dt] \\
&=E(u)(x)-[\alpha^{(4)}+M\alpha-N\int_0^1k(x,t)\alpha(t)dt] \\
&=[f(x,u,\int_0^1k(x,t)u(t)dt)+Mu-N\int_0^1k(x,t)u(t)dt]- [\alpha^{(4)}+M\alpha-N\int_0^1k(x,t)\alpha(t)dt]\\
&\geq M(u-\alpha)-N\int_0^1k(x,t)(u-\alpha)(t)dt+f(x,u,\int_0^1k(x,t)u(t)dt)-f(x,\alpha,\int_0^1k(x,t)\alpha(t)dt)\\
&\geq 0.\\
 \endaligned\eqno (3.8)
  $$
On the other hand, since $\Phi: C[0,1]\rightarrow D(L)$, then
$$
g(0)=(\Phi u)(0)-\alpha(0)=-\alpha(0)\geq 0, \ \ \ \ \ g(1)=(\Phi u)(1)-\alpha(1)=-\alpha(1)\geq 0 \eqno (3.9)
$$and
$$
g''(0)=(\Phi u)''(0)-\alpha''(0)=-\alpha''(0)\leq 0, \ \ \ \ \ g''(1)=(\Phi u)''(1)-\alpha''(1)=-\alpha''(1)\leq 0.\eqno (3.10)
$$ By maximum priciple in Theorem 2.2.\ (i),  (3.8)-(3.10) imply that $g\geq0$ and then $\alpha\leq \Phi(u)$.

 By a similar way, using the definition of the upper solutions, the maximum priciple in Theorem 2.2.\ (ii) and (3.4), we can get that $\Phi(u)\leq\beta$, then (3.7) is proved.

Based upon Schauder fixed-point theorem, $\Phi$ has a fixed point in $[\alpha,\beta]$ which is a solution of (1.1),(1.2).

Secondly, we show the following claim:  $$\beta\geq u_{1}\geq u_{2}\geq\alpha\Rightarrow \Phi u_{1}\geq\Phi u_{2}.\eqno (3.11)$$
In fact, let $\phi=\Phi u_{1}-\Phi u_{2}$, by (3.4) again, we have
$$
\aligned
&\phi^{(4)}+M\phi-N\int_0^1k(x,t)\phi(t)dt\\
&=[(\Phi u_{1})^{(4)}+M\Phi u_{1}-N\int_0^1k(x,t)(\Phi u_{1})(t)dt]-[(\Phi u_{2})^{(4)}+M\Phi u_{2}-N\int_0^1k(x,t)(\Phi u_{2})(t)dt] \\
&=E(u_{1})(x)-E(u_{2})(x)=[f(x,u_{1},\int_0^1k(x,t)u_{1}(t)dt)+Mu_{1}-N\int_0^1k(x,t)u_{1}(t)dt]\\
&\ \ \ \ \ \ \ \ \ \ \ \ \ \ \ \ \ \ \ \ \ \ \ \ \ \ \ \ \ \ \ \ \ \ \ -[f(x,u_{2},\int_0^1k(x,t)u_{2}(t)dt)+Mu_{2}-N\int_0^1k(x,t)u_{2}(t)dt]\\
&\geq 0.\\
 \endaligned\eqno (3.12)
  $$
On the other hand, $$
\phi(0)=\phi(1)=\phi''(0)=\phi''(1)=0. \eqno (3.13)
$$
Then by Remark 2.2., we conclude that $\phi\geq0$ and then the claim (3.11) is proved.

By the definition of the mapping $\Phi$, the iterative procedure (3.5) is equivalent to the iterative equation
$$
y_{n}=\Phi y_{n-1},\ \ \ \ n=1,2,\ldots.
$$
Define the iterative sequences $\{\alpha_{n}\}$ and $\{\beta_{n}\}$ satisfy $$
\alpha_{n}=\Phi\alpha_{n-1},\ \ \ \beta_{n}=\Phi\beta_{n-1},\ \ \ \ \ n=1,2,\ldots \eqno (3.14)
$$
with $\alpha_{0}=\alpha$ and $\beta_{0}=\beta$. Then combining (3.7) with (3.11), it is easy to see that $\{\alpha_{n}\}$ and $\{\beta_{n}\}$ have the monotonicity (3.6). By the compactness of $\Phi$ and the monotonicity (3.6), it follows that $\{\alpha_{n}\}$ and $\{\beta_{n}\}$ are convergent in $C[0,1]$, that is, there exist $\underline{y}$ and $\overline{y}\in C[0,1]$ such that $$\lim\limits_{n\rightarrow \infty}\alpha_{n}(x)=\underline{y}(x), \ \ \ \ \  \ \ \ \ \lim\limits_{n\rightarrow \infty}\beta_{n}(x)=\overline{y}(x).$$
On the other hand, it is easy to see that the operator $\Phi$ is continuous, then letting $n\rightarrow\infty$ in (3.14), we have
$$
\underline{y}=\Phi(\underline{y}),\ \ \ \ \ \ \ \ \  \overline{y}=\Phi(\overline{y}),
$$
thus $\underline{y}$ and $\overline{y}$ are the solutions of (1.1),(1.2).

Finally, we show that $\underline{y}$ and $\overline{y}$ are the extremal solutions of (1.1),(1.2) on $[\alpha, \beta].$

Let $y\in [\alpha,\beta]$ be an arbitrary solution of problem (1.1),(1.2), then combining (3.7) with (3.11) we have
 $$
\Phi^{n}\alpha\leq \Phi^{n}y\leq\Phi^{n}\beta,
$$
that is
$$
\alpha_{n}\leq y\leq\beta_{n}.
$$
Letting $n\rightarrow\infty$, we have
$$
\underline{y}\leq y\leq\overline{y}.
$$
Hence, $\underline{y}$ and $\overline{y}$ are minimum and maximum solutions of (1.1),(1.2) in $[\alpha, \beta]$ respectively.\hfill{$\Box$}

\section{Examples}

We present two examples to illustrate the application of Theorem 3.1.

\noindent{\bf Example 4.1.} Consider the nonlinear fourth order boundary value problem
$$
\aligned
&y^{(4)}(x)=\sin(\pi x)[(2-y^2(x))\int_0^1t y(t)dt+1],\ \ \ \ x\in (0,1),\\
&y(0) = y(1) = y''(0) = y''(1) = 0.\\
\endaligned
\eqno (4.1)
$$
Denote $$k(x,t)=\sin(\pi x)t, \ \ \ \ \ \ (x,t)\in [0,1]\times[0,1],$$ $$f(x,y(x),\int_0^1k(x,t)y(t)dt)=(2-y^2(x))\int_0^1\sin(\pi x)t y(t)dt+\sin(\pi x),\ \ \ \ x\in [0,1],$$
then $k\geq 0, \|k\|_{\infty}=1$ and the function $$f(x,u,v)=(2-u^2)v+\sin(\pi x)$$ is continuous on $[0,1]\times {\mathbb{R}}^{2}$. Moreover, problem (4.1) is equivalent to $$
\aligned
&y^{(4)}(x)=f(x,y(x),\int_0^1k(x,t)y(t)dt),\ \ \ \ x\in (0,1),\\
&y(0) = y(1) = y''(0) = y''(1) = 0.\\
\endaligned
\eqno (4.2)
$$
It is easy to see that $\alpha=0$ and $\beta=\sin(\pi x)$ are lower and upper solutions of (4.2) respectively which satisfy (3.3).

Taking $M=\frac{2}{\pi}$, by (2.3) and simple computation we have
$$\underset{x\in [0,1]}\max\int^{1}_{0}G(x,s)ds=\frac{\cosh \sqrt[4]{\frac{2}{\pi}}}{\frac{2}{\pi}},$$
then for any $0<N<\frac{\frac{2}{\pi}}{\cosh \sqrt[4]{\frac{2}{\pi}}}\approx 0.446$, it is easy to verify that (3.4) and other conditions in Theorem 3.1 are satisfied.

Hence by theorem 3.1, (4.2) has extremal solutions which can be obtained by iterative procedure (3.6). Moreover, since $y\equiv 0$ is not a solution of (4.1),  then the extremal solutions are positive.

\noindent{\bf Example 4.2.} Consider the nonlinear sixth order boundary value problem
$$
\aligned
&u^{(6)}(x)=f(x,u''(x),u(x)),\ \ \ \ x\in (0,1),\\
&u(0) = u(1) = u''(0) = u''(1) =u''''(0)=u''''(1)=0\\
\endaligned
\eqno (4.3)
$$where $f\in C([0,1]\times {\mathbb{R}}^{2}, \mathbb{R})$.
Let $u''(x)=y(x)$, then (4.3) is equivalent to (4.2) in which $$
k(x,t) =\left\{ \aligned
x(1-t), \ \ \ \ \ \ \ \  \ 0\leq x\leq t\leq1,\\
t(1-x), \ \ \ \ \ \ \ \  \ 0\leq t\leq x\leq1
\endaligned
\right.
$$
is the Green function for
$$\aligned
&u''(x)=0,\ \ \ \ x\in (0,1),\\
&u(0) = u(1) = 0.\\
\endaligned
$$
Obviously, $k\in C([0,1]\times [0,1], \mathbb{R})$ and $k\geq 0, \|k\|_{\infty}=\frac{1}{4}$, then problem (4.3) can be solved by the monotone iterative technique in Theorem 3.1 under certain conditions.

\vskip 3mm


\noindent{\bf Acknowledgements}

\noindent This work was supported by the NSFC(No.11961043,\ No. 11801453).

\vskip 12mm

\centerline {\bf REFERENCES}\vskip5mm\baselineskip 0.45cm
\begin{description}
\baselineskip 15pt

\item{[1]}~ T. V. Karman, M. A. Biot, Mathematical Methods in Engineering: An
Introduction to the Mathematical Treatment of Engineering Problems. McGraw¨CHill,
New York, 1940.


\item{[2]}~ A. Pugsley, The Theory of Suspension Bridge. Edward Arnold Pub. Ltd, London, 1968. 

\item{[3]}~ C. Itzykson, J. B. Zube, ¡±Quantum Field Theory,¡± Dover Publications,
INC, Mineola, New York, 2005.
\item{[4]}~ J. F. Kelly, Transient ultrasonic fields in power law attenuation media, 
    2008.
\item{[5]}~ E. Aruchunan, Y. Wu, B. Wiwatanapataphee, et al. A new variant of arithmetic mean iterative method for fourth order integro-differential equations solution, 2015 IEEE Third International Conference on Artificial Intelligence, Modelling and Simulation. IEEE, 2015.
\item{[6]}~ A. A. Dascioglu, M. Sezer, A taylor polynomial approach for solving high-order linear Fredholm integro-differential equations in the most general form, Int. J. Comput. Math. 
    84(4) (2007), 527-539.
\item{[7]}~ M. Lakestania, M. Dehghanb, Numerical solution of fourth-order integro-differential equations using Chebyshev Cardinal Functions, Int. J. Comput. Math. 87 (2010), 1389-1394.
\item{[8]}~ F. Ghomanjani, A. V. Kamyad, A. Kiliman,  Bezier curves method for fourth-order integro-differential equations, Abstr. Appl. Anal. 2013, 1-5.
\item{[9]}~ N. H. Sweilam, Fourth order integro-differential equations using variational iteration method, Comput. Math. Appl. 54 (2007),  1086-1091.
\item{[10]}~ A. Zeeshan, M. Atlas, Optimal solution of integro-differential equation of
suspension bridge model using genetic algorithm
and Nelder-Mead method, Journal of the Association of Arab Universities for Basic and Applied Sciences. 24 (2017), 310-314.
\item{[11]}~ W. M. Wang, An algorithm for solving the high-order nonlinear Volterra-Fredholm integro-differential equation with mechanization, Appl. Math. Comput. 172(1) (2006), 1-23.
\item{[12]}~ H. Amann, Fixed point equations and nonlinear eigenvalue problems in ordered Banach spaces.  SIAM Rev. 18 (1976), 620-709.
\item{[13]}~ G. S. Ladde, V. Lakshmikantham, A. S. Vatsala, Monotone Iterative Techniques for Nonlinear Differential Equations, Pitman, Boston, 1985.
\item{[14]}~ C. V. Pao, Nonlinear Parabolic and Elliptic Equations.  Plenum Press, New York, 1992.
\item{[15]}~ J. J. Nieto, An abstract monotone iterative technique, Nonlinear Anal. 28 (1997), 1923-1933.
    \item{[16]}~ I. Rachunkova, Upper and lower solutions and
topological degree, J. Math. Anal. Appl.  234 (1999), 311-327.
\item{[17]}~ C. D. Coster, P. Habets, The lower and upper solutions method for boundary value problems, in: A. Canada, P. Drabek and A. Fonda (Eds.), Handbook of Differential Equations-Ordinary Differential Equations, 2004.
\item{[18]}~ R. Agarwal, On fourth-order boundary value problems arising in beam analysis, Differ. Integr. Equ. 2 (1989), 91-110.

\item{[19]}~ C. De Coster, L. Sanchez, Upper and lower solutions, Ambrosetti-Prodi problem and positive solutions for fourth O.D.E, Riv. Mat. Pura Appl. 14 (1994), 1129-1138.

\item{[20]}~ A. Cabada, J. A. Cid, L. Sanchez, Positivity and lower and upper solutions for fourth order boundary value problems, Nonlinear Anal. 67 (2007), 1599-1612.

\item{[21]}~ R. Y. Ma, J. H. Zhang, S. M. Fu, The method of lower and upper solutions for fourth-order two-point boundary value problems, J. Math. Anal. Appl. 215 (1997), 415-422.

\item{[22]}~ R. Y. Ma, J. X. Wang, Y. Long, Lower and upper solution method for the problem of elastic beam, J. Fixed Point Theory Appl. 20 (2018) 1-13.
\item{[23]}~ R. Y. Ma, J. X. Wang, D. L. Yan, The method of lower and upper solutions for fourth order equations with the Navier condition, Bound. Value Probl. 152 (2017), 1-9.
\item{[24]}~ Y. X. Li, A monotone iterative technique for solving the bending elastic beam equations. Appl. Math. Comput. 217 (2010), 2200-2208.

\item{[25]}~D. X. Ma, X. Z. Yang, Upper and lower solution method for fourth-order four-point
boundary value problems, J. Comput. Appl. Math. 223 (2009), 543-551.

\item{[26]}~ M. H. Pei, S. K. Chang, Monotone iterative technique and symmetric positive solutions for a fourth-order boundary value problem, Math. Comput. Modelling. 51 (2010), 1260-1267.
\item{[27]}~ Y. M. Wang, Monotone iterative technique for numerical solutions of fourth-order
nonlinear elliptic boundary value problems, Appl. Numer. Math. 57 (2007), 1081-1096.

\item{[28]}~ Z. B. Bai, The method of lower and upper solutions for a bending of an elastic beam equation, J. Math. Anal. Appl. 248
(2000), 195-202.

\item{[29]}~ Q. L. Yao, On the positive solutions of Lidstone boundary value problems, Appl. Math. Comput. 137 (2003), 477-485.

\item{[30]}~ D. Q. Jiang, W. J. Gao, A. Y. Wan, A monotone method for constructing extremal solutions to fourth-order periodic
boundary value problems, Appl. Math. Comput. 132 (2002), 411-421.
\item{[31]}~ Q. Zhang, S. H. Chen, J. H. Lu,  Upper and lower solution method for fourth-order four-point boundary value problems, J. Comput. Appl. Math. 196 (2006), 387-393.

\item{[32]}~ J. Ehme, P. W. Eloe, J. Henderson,   Upper and lower solution methods for fully nonlinear
boundary value problems, J. Differ. Equations. 180 (2002), 51-64.

\item{[33]}~ P. Habets, L. Sanchez, A monotone method for fourth order boundary value problems involving a factorizable linear
operator, Port. Math. 64(3) (2007), 255-279.

\item{[34]}~ J. F. Fialho, F. Minhos, Higher order functional boundary value problems without monotone assumptions, Bound. Value Probl. (2013), 81.

\item{[35]}~ R. Vrabel, On the lower and upper solutions method for the problem of elastic beam with hinged ends. J. Math. Anal. Appl. 421(2) (2015), 1455-1468.

\item{[36]}~ E. Alves, T. F. Ma, L. P. Mauricio, Monotone positive solutions for a fourth order equation with nonlinear boundary conditions, Elsevier Science Publishers B. V. 2009.

\item{[37]}~ Q. A. Dang , T. K. Q. Ngo, Existence results and iterative method for solving the cantilever beam equation with fully nonlinear term, Nonlinear Anal. RWA. 36 (2017), 56-68.

 \item{[38]} J. A. Cid, D. Franco, F. Minhos, Positive fixed points and fourth-order equations, Bull. Lond. Math. Soc. 41 (2009), 72-78.

\end{description}
\end{document}